\documentclass[12pt]{article}
\usepackage{graphics}

\date{}
\def\qed{{\hfill\rule{1.2ex}{1.2ex}}}

\def\C{{\bf C}}
\def\R{{\bf R}}
\def\Q{{\bf Q}}

\def\Z{{\bf Z}}

\author{Norbert A'Campo}
\title{TQFT computations and experiments.}

\begin{document}

\maketitle

Our recent computer program TQFT 
allows the actual computation of the 
projective representation of the mapping class group of 
a surface on its Verlinde modules.
The program is build upon  the fusion formulae  of 
Roberts, Masbaum and Vogel [L],[M-V],[R].  The author
is greatly indebted to Gregor Masbaum for his 
explanations during the writing of
the program. The program uses the powerful calculator Pari-gp.

The actual program is mainly written for surfaces $\Sigma_{g,1}$
of genus $g\geq 1$ with one boundary component. 
We think of the surface $\Sigma_{g,1}$ as the boundary of a thickening 
$H_g$ in $\R^3$ 
of the trivalent graph $\Gamma_g$ of Fig. $1$. The graph $\Gamma_g$ has 
$3g-1$ edges. The Verlinde module $V^g_{k,i}$ appears as a vector space
over the field $C_{2k+4}=\Q(A)/(\Phi_{2k+4}(A))$, where $\Phi_{2k+4}$ is the
cyclotomic polynomial, whose roots are the primitive $(2k+4)$th roots of unity.
The space $V^g_{k,i}$ can be viewed as the $C_{2k+4}$ span 
of the set of admissible 
$(k,i)$-colorings of the edges of the graph $\Gamma_g$ by integers. 
An admissible $(k,i)$-coloring of $\Gamma_g$ is 
by definition an edge  coloring  $c:{\rm Edge}(\Gamma_g) \to \Z$
satisfying the following conditions:\newline
1. For each edge $e$ the integer $c(e)$ is even and 
satisfies $0\leq c(e) \leq k$;\newline
2. For each node of $\Gamma_g$ with adjacent edges $e_1,e_2,e_3$ the
colors $c(e_1),c(e_2),c(e_3)$ satisfy the triangular inequalities and the
inequality $c(e_1)+c(e_2)+c(e_3)\leq 2k$;\newline
3. The color of the outgoing edge is $i$.\newline

We describe the action of the mapping class group 
${\rm Mod}_{g,1}$ on $V^g_{k,i}$  
in terms of the actions of the following generators of ${\rm Mod}_{g,1}$.
Let $A_e$ be (up to isotopy) the simple loop on $\Sigma_{g,1}$ bounding in 
$H_g$ an embedded disk having one transversal intersection with edge $e$.
Let $B_r$ be the simple loop surrounding the handle $r$. 
The group ${\rm Mod}_{g,1}$
is generated by the Dehn twists about the $A_e$ and $B_r$. 

We explain some basic possibilities of the program.
See the 00README of {\tt http://www.geometrie.ch/TQFT} and also have  a
look at the explanations in the
comments of the pari-script tqft.gp.
The program TQFT  can only be called from a Pari session.

\begin{center}
\scalebox{1}[1]{\includegraphics{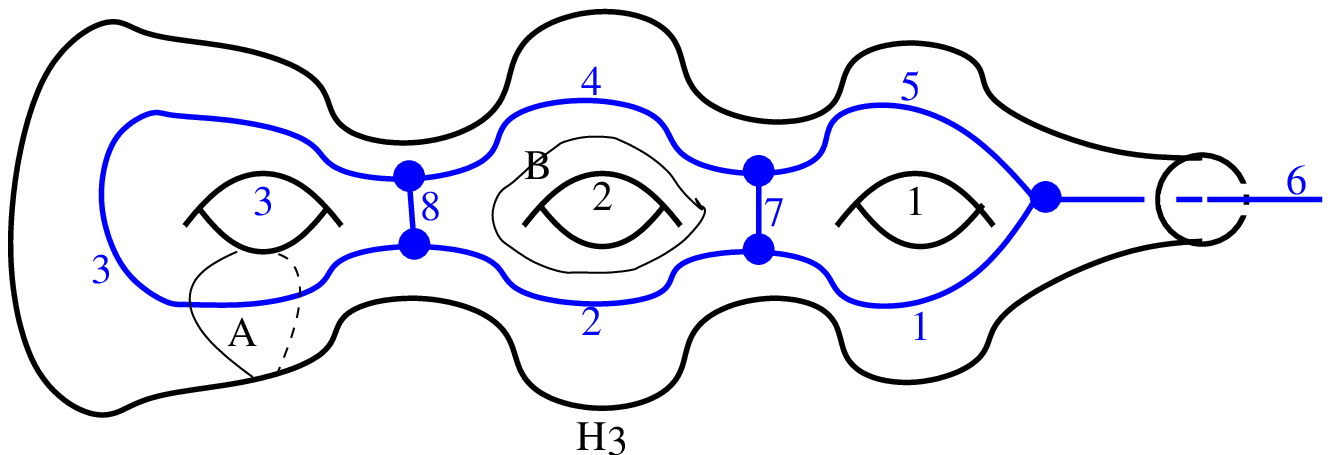}}
\newline
{Fig. 1. Handle body of genus 3. Trivalent graph with input edge.}
\end{center}

Now start a pari session 
Pari-gp with the command ``gp'' and read the file tqft.gp into
the pari session. 
The command

$init\_ so(k)$

\noindent
initializes the global variable $A$ and the field 
$C_{2k+4}=\Q(A)/(\Phi_{2k+4}(A))$.
The command

$init\_ boom\_ so([0,1, .. ,g-1],i)$

\noindent
initializes the graph $\Gamma_g$ and computes the list of its admissible
$(k,i)$-colorings. 
The commands 

$twA(e)$\newline\indent
$twB(r)$\newline\indent

\noindent
compute  matrices with coefficients in the field $C_{2k+4}$ 
of the action of the right Dehn twist about $A_e$ or 
$B_r$ in the space  $V_{k,i}$ with as basis 
the list of admissible colorings.

Our computations lead to  the following results:

{\bf Case genus $g=1$.} 

In this case the group ${\rm Mod}_{g,1}$ is identified via
its action on $H_1(\Sigma,\Z)$ with ${\rm SL}(2,\Z)$. The space $V^1_{k,i}$
is generated by the admissible colorings $(j,i)$ of the following graph

The admissibility conditions restrict to:
$$i \leq 2j, i+2j \leq 2k, \,j\, {\rm even}$$
hence we have the colorings
$$(j,i), j=i/2,i/2+2, \cdots , (k-1-i)/2$$
if $i$ and $k-1$ are divisible by $4$,
$$(j,i), j=i/2+1,i/2+3, \cdots , (k-1-i)/2$$,
if $i$ nor $k-1$ are divisible by $4$,
$$(j,i), j=i/2,i/2+2, \cdots , (k-3-i)/2$$
if $i$ divisible by $4$ but not $k-1$ and finally
$$(j,i), j=i/2+1,i/2+3, \cdots , (k-1-i)/2$$
if $k$ divisible by $4$ but not $i$.
So for $k=2h+1\geq 3$ the dimension of the vector space 
$V^1_{2h+1,2}$ is $h$. 

\begin{center}
\scalebox{1}[1]{\includegraphics{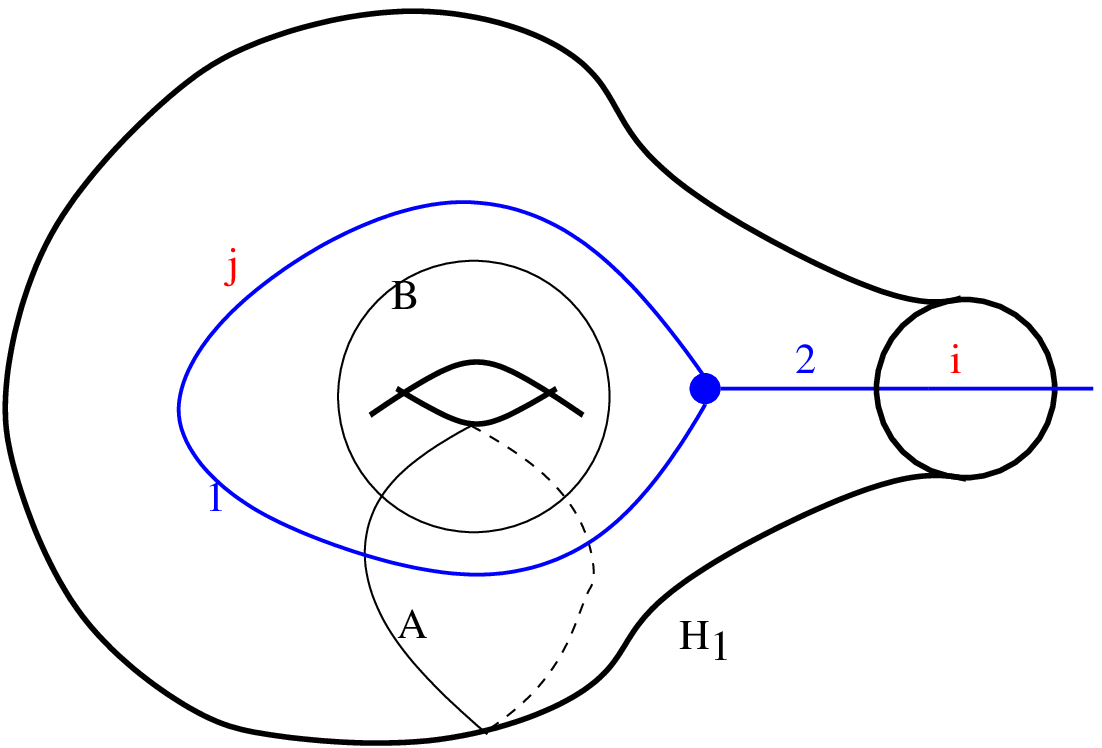}}
\newline
{Fig. $2$. Handle body of genus $1$. Loops $A$ and $B$.}
\end{center}

Let $D_1,D_2$ be the right Dehn twists about the curves 
$A_1,B_1$ of the Fig. $2$. The actions in homology $D_{1;*}$ and 
$D_{2;*}$ for a suitable 
orientation of the the torus are given by the matrices
$D_{1;*}=[1,1;0,1]$ and $D_{2;*}=[1,0;-1,1]$. The actions of  $V^1_{5,2}$
are obtained by doing first the commands:

$init\_ so(5);init\_ boom\_ so([0],2);$

\noindent
followed by 

$twA(1)$

\noindent
for $D_{1;5,2}$ and by

$twB(1)$

\noindent
for $D_{2;5,2}$.
One gets the matrices 
(remember that the coefficients are in $C_{14}$):

$$
D_{1;5,2}={\pmatrix{
 {A^{5} - A^{4} + A^{3} - A^{2} + A - 1}&{0}\cr
  {0}&{A^{4}}\cr
   }}
$$
$D_{1;5,2}=$
$$
{\pmatrix{
 {{{3}\over{7}}A^{5} + {{1}\over{7}}A^{4} + {{2}\over{7}}A^{3} - {{5}\over{7}}A^{2} + {{1}\over{7}}A - {{4}\over{7}}}&{{{-4}\over{7}}A^{5} + {{8}\over{7}}A^{4} - {{5}\over{7}}A^{3} + {{2}\over{7}}A^{2} - {{6}\over{7}}A + {{3}\over{7}}}\cr
  {{{-3}\over{7}}A^{5} + {{6}\over{7}}A^{4} - {{2}\over{7}}A^{3} + {{5}\over{7}}A^{2} - {{1}\over{7}}A + {{4}\over{7}}}&{{{4}\over{7}}A^{5} - {{1}\over{7}}A^{4} + {{5}\over{7}}A^{3} - {{2}\over{7}}A^{2} + {{6}\over{7}}A - {{3}\over{7}}}\cr
   }}
$$   

We now change the level to $k=7$ and work with the representation
$\rho_7:{\rm SL}(2,\Z) \to {\rm PGL}(V^1_{7,2})={\rm PGL}(3,C_{18})$.
This requires the commands:

$init\_ so(7);init\_ boom\_ so([0],2];$

The  image of the right Dehn twist about $A(1)$ is represented by 
the following matrix with coefficients in $C_{18}=\Q(A)/(A^6-A^3+1)$\newline 
$D_{1;7,2}=$
$$
{\pmatrix{
 {-A}&{0}&{0}\cr
  {0}&{-A^{3}}&{0}\cr
   {0}&{0}&{A^{3} - 1}\cr
    }}
$$

It is interesting to observe that the matrices $D_{1;7,2}$ and its
inverse $D_{1;7,2}^{-1}$
do not represent conjugate
elements in the group ${\rm PGL}(3,C_{18})$ , since the class function
$m\in {\rm PGL}(3,C_{18}) \mapsto {\rm trace}(m)^3/{\rm det}(m) \in C_{18}$
takes on  $D_{1;7,2}$ the value 
$$Mod(2A^5 - A^2, A^6 - A^3 + 1)$$ 
and 
on $D_{1;7,2}^{-1}$ the value 
$$Mod(-A^4 - A, A^6 - A^3 + 1).$$ 
It follows, as noticed by 
Vladimir Turaev, that the representation
$$
\rho_7:{\rm SL}(2,\Z) \to {\rm PGL}(V^1_{7,2})
$$ 
does not extend to a
representation on ${\rm GL}(2,\Z)$ since $[1,1;0,1]$ and $[1,-1;0,1]$
are in ${\rm GL}(2,\Z)$ conjugate.

Our next observation is that the images of the matrices 
$a=[7,3;2,1]$ and $b=[7,1;6,1]$ under 
$\rho_5:{\rm SL}(2,\Z) \to {\rm PGL}(V^1_{5,2})={\rm PGL}(2,\Q(A)/(\Phi_{14}(A))$ 
are not conjugate in ${\rm PGL}(2,\Q(A)/(\Phi_{14}(A))$. This is worth noticing
since the matrices $a,b$ are conjugate in ${\rm SL}(2,\Q)$.

From this observation we speculate about a positive answer to the following
question: given two elements $a,b \in {\rm SL}(2,\Z)$ that 
are not conjugate in 
${\rm SL}(2,\Z)$ does there exist $k\geq 3, k\, {\rm odd}$, such that
the images $\rho_k(a),\rho_k(b)$ are not conjugate in 
${\rm PGL}(2,\Q(A)/(\Phi_{2k+4}(A))$?

We wish to compute for $a=[7,3;2,1]$ and $b=[7,1;6,1]$ in level 
$k=5$ and with input color $i=2$, so 
initialize back to $k=5,i=2$ with:\newline
$init\_ so(7);init\_ boom\_ so([0],2];$\newline
The following commands create the matrices $a,b$:\newline
$a=Mat([7,3;2,1]);b=Mat([7,1;6,1]);$\newline
With the commands:\newline
$wa=slw(a);wb=slw(b);$\newline
we write $a,b$ as products of the matrices 
$D_{1;*}=[1,1;0,1]$ and $D_{2;*}=[1,0;1,-1]$. 

$$wa=
{\pmatrix{ {\pmatrix{
 {1}&{1}\cr
  {0}&{1}\cr
   }}&{\pmatrix{
    {1}&{1}\cr
     {0}&{1}\cr
      }}&{\pmatrix{
       {1}&{1}\cr
        {0}&{1}\cr
	 }}&{\pmatrix{
	  {1}&{0}\cr
	   {1}&{1}\cr
	    }}&{\pmatrix{
	     {1}&{0}\cr
	      {1}&{1}\cr
	       }}\cr}
	       }
$$
$$
wb={\pmatrix{ {\pmatrix{
 {1}&{1}\cr
  {0}&{1}\cr
   }}&{\pmatrix{
    {1}&{0}\cr
     {1}&{1}\cr
      }}&{\pmatrix{
       {1}&{0}\cr
        {1}&{1}\cr
	 }}&{\pmatrix{
	  {1}&{0}\cr
	   {1}&{1}\cr
	    }}&{\pmatrix{
	     {1}&{0}\cr
	      {1}&{1}\cr
	       }}&{\pmatrix{
	        {1}&{0}\cr
		 {1}&{1}\cr
		  }}&{\pmatrix{
		   {1}&{0}\cr
		    {1}&{1}\cr
		     }}\cr}
		     }
$$		     

The corresponding products of the matrices
$D_{1;5,2}$ and $D_{2;5,2}$ computes the actions $Va$ and $Vb$ 
of $a$ and $b$ on
$V^1_{1,2}$.
On gets these matrices directly with the following commands:

$Va=eval\_ sl(a,5,2);Vb=eval\_ sl(b,5,2);$

$Va=$
$$
{\pmatrix{
 {{{8}\over{7}}A^{5} - {{2}\over{7}}A^{4} + {{3}\over{7}}A^{3} - {{4}\over{7}}A^{2} + {{5}\over{7}}A - {{6}\over{7}}}&{{{1}\over{7}}A^{5} + {{5}\over{7}}A^{4} - {{4}\over{7}}A^{3} + {{3}\over{7}}A^{2} - {{2}\over{7}}A + {{1}\over{7}}}\cr
  {{{6}\over{7}}A^{5} + {{2}\over{7}}A^{4} + {{4}\over{7}}A^{3} - {{3}\over{7}}A^{2} + {{2}\over{7}}A - {{8}\over{7}}}&{{{6}\over{7}}A^{5} + {{2}\over{7}}A^{4} + {{4}\over{7}}A^{3} - {{3}\over{7}}A^{2} + {{2}\over{7}}A - {{1}\over{7}}}\cr
   }}
$$

$Vb=$
$$
{\pmatrix{
 {{{-4}\over{7}}A^{5} + {{1}\over{7}}A^{4} - {{5}\over{7}}A^{3} + {{2}\over{7}}A^{2} + {{1}\over{7}}A + {{3}\over{7}}}&{{{3}\over{7}}A^{5} + {{1}\over{7}}A^{4} - {{5}\over{7}}A^{3} + {{2}\over{7}}A^{2} + {{1}\over{7}}A + {{3}\over{7}}}\cr
  {{{4}\over{7}}A^{5} - {{1}\over{7}}A^{4} + {{5}\over{7}}A^{3} - {{2}\over{7}}A^{2} - {{1}\over{7}}A - {{3}\over{7}}}&{{{4}\over{7}}A^{5} - {{1}\over{7}}A^{4} - {{2}\over{7}}A^{3} - {{2}\over{7}}A^{2} - {{1}\over{7}}A - {{3}\over{7}}}\cr
   }}
$$

Our claim is that $Va,Vb$ are not conjugate in ${\rm PGL}(2,C_{14})$:
The quantities 
$${{trace(Va)^2} \over {det(Va)}}\in C_{14},\,\, 
{{trace(Vb)^2} \over {det(Vb)}}\in C_{14}
$$ 
are conjugacy invariants in 
the group ${\rm PGL}(2,C_{14})$
and we deduce our the claim from 

${{trace(Va)^2} \over {det(Va)}}=
Mod(-3A^5 + 2A^4 - 2A^3 + 3A^2 + 5, \Phi_{14}(A)),$

${{trace(Vb)^2} \over {det(Vb)}}=Mod( 1,\Phi_{14}(A)).$

The representations 
$\rho_k:{\rm SL}(2,\Z) \to {\rm PGL}(2,\Q(A)/(\Phi_{2k+4}(A)))$ 
can be lifted to ${\rm GL_1}(2,\Q(A)/(\Phi_{2k+4}(A))/\mu_{2k+4}$,
hence the maximal absolut value of  $\sigma({\rm trace}(\rho_k(a)))$,
$\sigma$ running over all the embeddings of the field 
$\Q(A)/(\Phi_{2k+4}(A))$ into $\C$, is an invariant $||a||_k$ 
for $a$. Here we have denoted by ${\rm GL_1}$ the subgroup of elements
in ${\rm GL}$ having a root of unity as determinant.

The maximal absolut value of  $\sigma({\rm trace}(\rho_5(a)))$ and of
$\sigma({\rm trace}(\rho_5(b)))$ are computed with the commands:

$normk(trace(Va));normk(trace(Vb));$

We get\newline 
$||a||_5=2.8019377358048382524722046390148901023$\newline
$||b||_5=1$,\newline 
showing once more that
$Va,Vb$ are not conjugate in ${\rm PGL}(2,C_{14})$.

Our second experiment concerns the growth of the TQFT action of 
the matrix $a=[2,1;1,1]$.
We have observed that 
the inequalities $||a||_k < trace(a)=3,k=3,7,9, \cdots $
hold,
and that $3$ is the supremum of $\{||a||_k\ \mid k\,\, {\rm odd}\}$.
Here the output of the following command line:

\noindent
$(15:12)\, gp\, >\, for(i=1,16,k=2*i+1;$\newline
$print(normk(trace(eval\_ sl([2,1;1,1],k,2))) ));$

\noindent
$1$\newline
$2.2469796037174670610500097680084796213$\newline
$2.8793852415718167681082185546494629398$\newline
$2.9189859472289947797807361141326553981$\newline
$2.7709120513064197918007510440301977572$\newline
$1$\newline
$2.8649444588087116091462317836431267725$\newline
$2.9727226068054447472050183896381342215$\newline
$2.9776616524502570901394857658680172261$\newline
$2.9258345746955985900304471947464775986$\newline
$1$\newline
$2.9460897411596476776657703455693918400$\newline
$2.9882759143087192179106054317591031337$\newline
$2.9897386467837902926427066197674389859$\newline
$2.9638573945254134007973488852494919219$\newline
$1$\newline
$time = 23mn,\, 6,060\, ms.$

Since $V_{5,2}$ has dimension $2$, we may deduce from
the second line of output \newline
$||a||_5=2.246...$, 
that the action of
$a=[2,1;1,1]$ on $V^1_{5,2}$ is a very simple and explicit 
example of an element of  
infinite order in TQFT, see [F],[Gi],[M]. Indeed,
the product of the two eigenvalues $\lambda_1,\lambda_2$
of $\rho_5(a)$ is a $14$th root of unity and for some embedding
$\sigma$ of the field that contains $\lambda_1,\lambda_2$
we have $|\sigma(\lambda_1)+\sigma(\lambda_2)|>2$;
it follows ${\rm max}(|\sigma(\lambda_1)|,|\sigma(\lambda_2)|)>1$.
It follows from the corollary of Theorem $2$ [Gi] of P. Gilmer
that the action of $a$ on $V^1_{k,0}$ any $k$  are periodic.

{\bf Case $g>1$.} 

The two slalom knots $K_1$ and $K_2$, 
see [AC1,AC2], of the 
planar rooted trees $[0,1,1,1,2]$ and $[0,1,1,1,3]$
in the table of KNOTSCAPE [H-T] 
are the knots $13n1320$ and $13n1291$ respectively.
According to the author's experience, it is 
difficult  to separate this pair of mutant knots by invariants.
They have for instance, equal  Kauffman polynomial and HOMFLY polynomial.
Moreover, the Khovanov homologies coincide, as one can verify using
the program KhoHo of Alexander Shumakovitch, see
{\tt http://www.geometrie.ch/KhoHo}. With SNAPPEA [W] we were able 
to show that these knots have non isomorphic rigid symmetry groups. 
With SNAP [G] we could not find a distinction based on arithmetic 
properties of the
hyperbolic structures on the complements. The knots are
fibered with fibers of genus $5$ and 
the monodromy diffeomorphisms can be
written explicitly  in terms of the underlying planar rooted trees 
as products of Dehn twist. 

The program TQFT was written and especially designed in order 
to compute the 
action of monodromies of slalom knots in TQFT. The following 
command line computes the TQFT actions for level $k=3$ and input color $i=2$
on $V^5_{3,2}$:
\newline
$c_1=coxeter([0,1,1,1,2],3,2);c_2=coxeter([0,1,1,1,3],3,2);$\newline
The compution is finished after: $time =\, 5mn,\, 930\, ms$.
We do not ask for an output on the screen, since $c_1,c_2$ are 
square matrices of size $275$ representing elements in 
${\rm PGL}(275,C_{10})$. First we compute the traces and get:\newline
$trace(c_1)=Mod(7A^3 + A^2 + 5A - 4, A^4 - A^3 + A^2 - A + 1) $,\newline
$trace(c_2)=Mod(7A^3 + A^2 + 5A - 4, A^4 - A^3 + A^2 - A + 1)$.\newline
So we have equal non vanishing traces. We compute traces of iterates:\newline
$trace(c_1^2)=Mod(-13A^3 + 18A^2 + 4A + 25, A^4 - A^3 + A^2 - A + 1) $,\newline
$trace(c_2^2)=Mod(-13A^3 + 18A^2 + 4A + 25, A^4 - A^3 + A^2 - A + 1) $.\newline
Again equal, so, we continue with traces of third powers:\newline
$trace(c_1^3)=Mod(-47A^3 - 56A^2 - 65A - 2, A^4 - A^3 + A^2 - A + 1) $,\newline
$trace(c_2^3)=Mod(-62A^3 - 47A^2 - 68A - 6, A^4 - A^3 + A^2 - A + 1) $.\newline
Third power traces are different, so since the traces of first powers 
are equal and non-vanishing, 
we conclude that the TQFT-actions of the monodromies of 
the knots $K_1,K_2$ are not
conjugate in ${\rm PGL}(275,C_{10})$.

We wish to  compute for the two knots $K_1,K_2$ also with input color $i=0$,
so we  use the commands:\newline
$c_1=coxeter([0,1,1,1,2],3,0);c_2=coxeter([0,1,1,1,3],3,0);$\newline
which was done after  $time = 2mn,\, 18,840\, ms.$
This corresponds to a study of the monodromy as a diffeomorphism of the
closed surface of genus $5$. We now get  matrices of smaller size $175$
and as before, the  separation of 
the knots with traces of third powers:\newline
$trace(c_1)=Mod(5A^3 + 3A - 3, A^4 - A^3 + A^2 - A + 1) $,\newline
$trace(c_2)=Mod(5A^3 + 3A - 3, A^4 - A^3 + A^2 - A + 1)$,\newline
$trace(c_1^2)=Mod(-8A^3 + 11A^2 + 3A + 15, A^4 - A^3 + A^2 - A + 1) $,\newline
$trace(c_2^2)=Mod(-8A^3 + 11A^2 + 3A + 15, A^4 - A^3 + A^2 - A + 1) $,\newline
$trace(c_1^3)=Mod(-20A^3 - 37A^2 - 35A - 25, A^4 - A^3 + A^2 - A + 1) $,\newline
$trace(c_2^3)=Mod(-29A^3 - 31A^2 - 38A - 13, A^4 - A^3 + A^2 - A + 1) $.\newline

The two slalom knots of the planar rooted trees $[0,1,1,3,3,4]$ and\newline
$[0,1,1,3,3,5]$ are the knots $15n30444$ and $15n30419$ respectively.
For this pair we have again that 
the traces of the first and second powers 
of the actions of the monodromies in level $k=3$ with input color $i=0$
on $V^6_{3,0}$
coincide, and that the traces of the third powers
distinguish the knots of this pair. For this 
calculation the size of the matrices grew up to $675$.

\noindent
[AC1]
Norbert A'Campo,
{\it Generic immersions of curves, knots,
monodromy and gordian number},
Publ. Math. I.H.E.S. {\bf 88} (1998), 151-169, (1999).

\noindent
[AC2]
Norbert A'Campo,
{\it Planar trees, slalom curves and hyperbolic knots},
Publ. Math. I.H.E.S. {\bf 88} (1998), 171-180, (1999).

\noindent
[F]
Louis Funar,
{\it On TQFT representations of the mapping class groups},
Pacific J. of Math. {\bf 188} 1999, 251--274.

\noindent
[Gi]
Patrick S. Gilmer,
{\it On the Witten-Reshetikhin-Turaev representations of mapping
class groups}, {\tt http://arXiv.org/q-alg/9701042\_ v5}.

\noindent
[G]
Oliver Goodman, {\it Snap},\newline {\tt http://www.ms.unimelb.edu.au/\~ snap}

\noindent
[H-T],
Jim Hoste and Morven Twistlethwaite,
{\it KNOTSCAPE},\newline {\tt http://dowker.math.utk.edu/knotscape.html}

\noindent
[L]
W.B.R. Lickorish, 
{\it Skeins and handlebodies}, Pac. J. Math.
  {\bf 159} (1993), 337-350.

\noindent
[M]
Gregor Masbaum, 
{\it An element of infinite order in TQFT-representations
   of mapping class groups}, Cont. Math. {\bf 223} (1999), 137-139.

\noindent
[M-V]
Gregor Masbaum, Pierre Vogel,
{\it $3$-valent graphs and the Kauffman bracket}, 
Pacific J. Math.  {\bf 164}  (1994),  no. 2, 361--381.

\noindent
[R]
J.D. Roberts, 
{\it Skeins and mapping class groups}, Math. Proc.
  Camb. Phil. Soc. {\bf 115} (1994), 53-77.

\noindent
[S]
Alexander Shumakovitch, {\it KhoHo},\newline {\tt http://www.geometrie.ch/KhoHo}

\noindent
[W]
Jeff Weeks, {\it SNAPPEA}, \newline{\tt http://www.northnet.org/weeks}

\noindent
Universitaet Basel, Rheinsprung 21, CH-4051  Basel.\newline
{\tt email: Norbert.Acampo@unibas.ch}
\end{document}